\documentclass[a4paper,11pt]{article}
\usepackage{amsmath,amsfonts,latexsym}
\usepackage{amsthm}
\usepackage{amssymb}
\usepackage[dvips]{graphicx}
\usepackage{array}

\makeatletter

\newcommand{\C}{\mathbb C}
\newcommand{\R}{\mathbb R}
\newcommand{\Z}{\mathbb Z}
\newcommand{\N}{\mathbb N}

\theoremstyle{plain}
\newtheorem{Thm}{Theorem}

\newtheorem{Pro}[Thm]{Proposition}
\newtheorem{Cor}{Corollary}[Thm]
\newtheorem{Lem}[Thm]{Lemma}

\theoremstyle{remark}

\newtheorem{Rem}[Thm]{\sc Remark}

\newcommand{\dpcm }{{\hfill$\Box$}\medskip}

\newcommand{\ps@taiplain}{%
     \renewcommand{\@oddhead}{\centerline{\sf {\small On the simpleness of zeros of Stokes multipliers}}}%
     \renewcommand{\@evenhead}{\@oddhead}%
     \renewcommand{\@oddfoot}{\centerline{\thepage}}
     \renewcommand{\@evenfoot}{\@oddfoot}}
\makeatother     

\begin{document}

\thispagestyle{plain}

\title{\bf On the simpleness of zeros of Stokes multipliers}
\vspace{-2cm}
\author{
{\em TRINH Duc Tai}
\thanks{Department of Mathematics, Teacher Training College of Dalat,
  29 Yersin, Dalat, Vietnam / tel.
84 63 827344 / fax. 84 63 834732/  {\sc email:}
tductai@hcm.vnn.vn}'\thanks{The Abdus Salam International Centre
for Theoretical Physics, Strada Costiera 11, Trieste 34014, Italy/
{\sc email:} tductai@ictp.trieste.it}} \vspace{-2cm}
\date{\today}

\maketitle

\pagestyle{taiplain}

\begin{abstract}
The aim of this paper is to discuss the simpleness of zeros of
Stokes multipliers associated with the differential equation
$-\Phi''(X) + W(X)\Phi(X) =0$, where $ W(X) = X^{m}
+a_{1}X^{m-1}+\cdots +a_{m}$ is a real monic polynomial. We show
that, under a suitable hypothesis on the coefficients $a_k$, all the
zeros of the Stokes multipliers are simple.
\end{abstract}



\section{Introduction}
We consider in the complex plane a second-order linear
differential equation,
\begin{equation}\label{EquatX}
-\Phi''(X) + W(X)\Phi(X)=0,
\end{equation}
where $ W(X) = X^{m} +a_{1}X^{m-1}+\cdots +a_{m}$ is a monic
polynomial of degree $m \in \N$.

Equation (\ref{EquatX}) has an irregular singular point at $X =
\infty$, so that the asymptotic behaviors of the solutions at this
point usually exhibit Stokes phenomena which are controlled by the
Stokes multipliers, in relation with the so-called connection
formulae.

These Stokes multipliers measure the lack of commutativity between
the analytic continuations of the solutions and their asymptotics
near infinity, when crossing the Stokes singular directions. This
justified their systematic study which has been initiated by
Sibuya in his book \cite{Sibuya} and  extended  in many studies,
in particular in the framework of the resurgent asymptotic
analysis (see, e.g., \cite{DDP93,DDP97,DelabR,Ph96,Ph00,Tai1}).

Since equation (\ref{EquatX}) depends on the parameter $(a_1,
\cdots, a_m)$, this translates to the Stokes multipliers. For
a convenient normalization, the Stokes multipliers are in fact
holomorphic functions in  $(a_1, \cdots, a_m)$, and the question
of describing their zeros appears as a natural mathematical
question. Specializing this question in the sole parameter $a_m$,
 Sibuya has shown that all but
a finite number of the zeros are simple (see \cite{Sibuya}, ch.
6), by exploring the asymptotic expansion at infinity of the Stokes
multipliers with respect to  $a_m$.

In this text we would like to extend this result, considering the
question of the simpleness of {\em all}  these zeros. This question
arises from the fact that the  zeros of
the Stokes multiplier are nothing but eigenvalues of a (complex)
boundary value problem associated with equation (\ref{EquatX}). Such a
problem has merged recently in the context of the physically well
motivated study of the so-called $\mathcal{PT}$-symmetric models (see,
e.g., \cite{Ben1,BenBoe99}).

In the paper \cite{Tai3},  we have shown that the
 simpleness of all the zeros of a convenient Stokes multiplier implies
the non-degeneration of the eigenstates under the so-called
$\mathcal{PT}$-pseudo-norm. It is commonly believed that this
property allows to define a conventional mathematical structure
for a physically consistent $\mathcal{PT}$-symmetric quantum
mechanic theory.

This explains the motivation of our present paper, that is to give
a complete proof
 for the simpleness of all the zeros of the Stokes
multipliers, under an appropriate condition  (from \cite{Shin002}) on the parameter
$(a_1, \cdots, a_{m-1})$ so as to impose the reality of the zeros.

As a preparation for the necessary ingredients, the next section
will serve to recall some important notions and facts from the
theory of Sibuya for equation (\ref{EquatX}). The most crucial
result is the existence and uniqueness theorem for the solutions
characterized by an appropriate asymptotic behavior at infinity.
Section 3 provides a detailed proof, in the spirit of \cite{Shin002},
 for the simpleness of zeros in an instructive special case, where all
$a_k$ (except for the last coefficient $a_m$) are vanishing.
A similar result, which holds for some more general cases, is then established  by the same
arguments. Finally, in
the conclusion, we briefly discuss a model exhibiting non-simple
real zeros.

\section{Stokes multipliers} \label{dapressib} In this section,
we briefly recall some classical results of Sibuya's theory on
second-order linear differential equations with polynomial
coefficients \cite{Sibuya}. The most fundamental fact is the
following theorem, which asserts the existence and uniqueness of a
solution characterized by its asymptotic behavior at infinity.

\begin{Thm}[Sibuya]\label{ketquaSibuya}
 Equation (\ref{EquatX}) admits a unique solution
 $ \Phi_0(X,a) := \Phi_0(X,a_{1},a_{2},\ldots,a_{m}) $ such that:
\begin{enumerate}
\item $\Phi_0(X,a)$ is an entire function in
                  $(X,a_{1},a_{2},\ldots,a_{m})$,
\item $\Phi_0(X,a)$ and its derivative $ \Phi'_0(X,a) $ admit the
following asymptotic behaviors
\begin{equation}\label{AsymptXprime}
\Phi_0(X) \simeq X^{r_{m}} e^{-S(X,a)} \left[ 1+O(X^{-1/2})\right]
\end{equation}
\begin{equation}\label{tiemcandaoham}
\Phi'_0(X) \simeq X^{\frac{m}{2}+r_{m}} e^{-S(X,a)} \left[
  -1+O(X^{-1/2})\right]
\end{equation}
when $X \rightarrow \infty$ in each sub-sector strictly contained
in the sector
$$ \Sigma_{0} = \left\lbrace \vert \arg(X) \vert <
\frac{3\pi}{m+2} \right\rbrace
$$
and the asymptotic regimes occur uniformly with respect to $a =
(a_{1},a_{2},\ldots,a_{m})$ in any compact of $\C^m$.
\end{enumerate}
\end{Thm}

In the above theorem $r_{m}$ and $S(X,a)$ can be determined
explicitly from $W(X)$. More concretely, as $X \rightarrow
\infty$, one can write
\begin{equation}\label{tam}
\begin{array}{ll}
  \sqrt{W(X)} & = X^\frac{m}{2}\left\{1 + a_{1}X^{-1}+\cdots
  +a_{m}X^{-m}\right\}^{1/2} \\
   & = X^\frac{m}{2}\left\{1+ \sum_{k=1}^{\infty}
 b_{k}(a)X^{-k}\right\}\\
\end{array}
\end{equation}
where, obviously, $b_{k}(a)$ are quasi-homogeneous polynomials in
$a_1,\ldots,a_m$ with real coefficients.

By integrating term-by-term the series in the right-hand side, we
get a primitive of $\sqrt{W(X)}$. The function $S(X,a)$ is
associated to the ``principal part" of this primitive
$$S(X,a) =
\frac{2}{m+2}X^{\frac{m+2}{2}} + \cdots $$
 that only contains
terms with strictly positive powers of $X$. And $r_m = r_m(a)$ is
defined by
\begin{equation}\label{dncuarm}
    r_m(a) = \left\{\begin{array}{ll}
      -m/4 & \rm{for}\; $m$\; \rm{odd} \\
       -m/4 - b_{1+m/2}(a) & \rm{for}\; $m$\; \rm{even}\\
    \end{array}\right.
\end{equation}

We should notice that for $m
> 2$, $r_m(a)$ does not depend on the last
coefficient $a_m$ and if all $a_j$ (possibly except $a_m$) are
equal to zero then $r_m =-m/4$.

We shall define other solutions of (\ref{EquatX}) by introducing a
rotation of the complex plan. Let us denote
$$\omega := e^{\frac{i2\pi}{m+2}}\quad {\rm and}\quad
 \omega_k(a):= (\omega^{k} a_1, \omega^{2k} a_2, \ldots ,
\omega^{km} a_m);\quad (k\in \Z)$$
For each $k\in \Z$, we construct functions $\Phi_k(X,a)$ by
setting
\begin{equation}\label{lesPhik}
\Phi_{k}(X,a) := \Phi_0(\omega^{-k}X,\omega_{-k}(a)).
\end{equation}

It is not difficult to check that $\Phi_{k}(X,a)$ are indeed
solutions of (\ref{EquatX}) and exponentially vanishing at
infinity in the corresponding sector
\begin{equation}\label{cacsector}
    S_k := \Big\{\big|\arg(X)-\frac{k2\pi}{m+2}\big| < \frac{\pi}{m+2}\Big\}.
\end{equation}
The following lemma, which can be verified in a straightforward
way (see \cite{DelabR,Sibuya}), implies the linear independence of
two consecutive solutions $\Phi_k$ and $\Phi_{k+1}$.
\begin{Lem}\label{bode1}
 For any $k \in \mathbb{Z}$, the Wronskian of $\Phi_k$ and $\Phi_{k+1}$ is given by the formula
  \begin{equation}
    \label{Wronskien}
{\sf Wr}(\Phi_{k},\Phi_{k+1}) = 2(-1)^k \omega ^ {\frac{km}{2}-
r_m(\omega_{-k-1}(a))}
  \end{equation}
\end{Lem}

From this observation, together with classical results on the
structure of solutions of linear differential equations, we can
infer that $\{\Phi_k,\Phi_{k+1}\}$ constitutes a basis for the
space of  solutions of  equation (\ref{EquatX}). Therefore, every
 solution can be expressed as a linear combination of
$\Phi_k,\Phi_{k+1}$. In particular, for each $k\in \Z$, we have
\begin{equation}\label{bieudien}
    \Phi_{k-1} =C_{k}(a)\Phi_{k}+\widetilde{C}_{k}(a)\Phi_{k+1}
\end{equation}
The multipliers $C_{k}(a)$ and $\widetilde{C}_{k}(a)$ are called
the \emph{Stokes multipliers} of $\Phi_{k-1}$ with respect to
$\Phi_{k}$ and $\Phi_{k+1}$. Further studies on these objects are
addressed in \cite{DelabR,Ph96,Sibuya}. By definition, it is
evident that
$$C_{k}(a)=\frac{{\sf Wr}(\Phi_{k-1},\Phi_{k+1})}{{\sf Wr}(\Phi_{k},\Phi_{k+1})}
\quad {\rm and }\quad \widetilde{C}_{k}(a) =\frac{{\sf
Wr}(\Phi_{k-1},\Phi_{k})}{{\sf Wr}(\Phi_{k+1},\Phi_{k})}
$$
Since $\Phi_k(X,a)$ are entire functions, it follows immediately
from these equalities and Lemma~\ref{bode1} that $C_k(a)$ and
$\widetilde{C}_k(a)$ are also entire functions in $a$. By
definition, $C_k$ is closely related to $C_0$ in a ``cyclic" way
through the formula:
$$C_k(a) = C_0(\omega_{k}(a))\, ,$$
by which the information about $C_k$ can be derived from $C_0$.
Furthermore, we also get an explicit expression for
$\widetilde{C}_k(a)$
$$
    {\widetilde{C}}_{k}(a) =  \omega ^
{-m-2r_m(\omega_{-k}(a))}.
$$
We emphasize that $\widetilde{C}_k(a)$ is never vanishing. It thus
can be reduced to 1 by a suitable renormalisation of the
$\Phi_k$'s. To do this, it is sufficient to insert a simple factor
in the right hand side of (\ref{lesPhik}). For instance, when
$k=0$, by redefining
\begin{equation}\label{ttt}
\begin{array}{ll}
   & \Phi_1(X,a) := \omega ^ {-m/2-r_m(a)}\Phi_0(\omega^{-1} X,\omega_{-1}(a)) \\\vspace{2.5mm}
  {\rm and} & \Phi_{-1}(X,a) := \omega ^ {m/2 + r_m(a)}\Phi_0(\omega X,\omega(a))\\
\end{array}
\end{equation}
we can write (\ref{bieudien}) under a slightly symmetric form,
\begin{equation}\label{ourcasedx}
    \Phi_{-1} =C(a)\Phi_{0}+\Phi_{1}\, ,
\end{equation}
where $C(a) := \omega ^ {m/2+r_m(a)}C_0(a)$ is also called the
\emph{Stokes multiplier} of $\Phi_{-1}$ with respect to $\Phi_0$.

Concerning this (sole) Stokes multiplier $C(a)$, whose zeros are
expected to be simple, we first have:

\begin{Pro}\label{md2} For any $a\in \C^m$,
\begin{equation}\label{hethuc}
             \overline{C(a)} + C(\overline{a}) = 0.
\end{equation}
\end{Pro}

\begin{proof}
By virtue of the quasi-homogeneity of equation (\ref{EquatX}), we
can see that $\overline{\Phi_0(\overline{X},\overline{a})}$ is
also one of its solutions whose asymptotic behavior at infinity in
the sector $S_0$ is the same as that of $\Phi_0(X,a)$.

The uniqueness of the canonical solution in Theorem
\ref{ketquaSibuya} implies immediately that
\begin{equation}\label{doixung}
    \overline{\Phi_0(X,a)} = \Phi_0(\overline{X},\overline{a}).
\end{equation}
Taking into account the above new definitions of $\Phi_{\pm 1}$
 in (\ref{ttt}), we can check without difficulty that
\begin{equation}\label{dxcuanghiem}
   \Phi_{-1}(\overline{X},\overline{a}) = \overline{\Phi_1(X,a)}
\end{equation}
for any $X\in \C$ and any $a\in \C^m$.

Putting these relations in (\ref{ourcasedx}) leads the desired
identity.
\end{proof}
\begin{Cor}\label{}
    The zero set of $C(a)$ is invariant under the complex conjugation $a \mapsto
    \overline{a}$.
\end{Cor}
\begin{Cor}\label{}
    Restricted on real coordinates
    $a = (a_1,\dots,a_m) \in\R^{m}$, $-iC(a)$ is a real-valued function.
\end{Cor}

\section{The simpleness of zeros of Stokes multiplier}
In this section, we shall discuss the simpleness of zeros of
$C(a)$ considered as a function of the last coefficient $a_m$. For
convenience, we consider $\lambda := a_m$ as variable of the
entire function $C(a,\lambda) := C(a_1,\ldots,a_{m-1},\lambda)$.
We then show that under some hypotheses on $a_1,
a_2,\ldots,a_{m-1}$, the derivative
$\frac{\partial}{\partial\lambda}C(a,\lambda)\neq 0$ if
$C(a,\lambda)=0$.

\subsection{A special case: $a_1=0,\, \ldots,a_{m-1}=0$}
We first concentrate on the case where all $a_j$ are vanishing
except for $a_m =: \lambda$. Equation (\ref{EquatX}) now simply
reads
\begin{equation}\label{donthuc}
    -\Phi''(X) + (X^m + \lambda)\Phi(X)=0
\end{equation}
Even for this simple case, the study of the zeros of the Stokes
multiplier is interesting because it has an intimate relation with the
spectral analysis problem of  Hamiltonians whose potentials are (possibly complex)
homogeneous polynomials (see
\cite{BenBoe98,BenBoe98-2,Vo83,Vor99-2}).

We shall use all the notations of the previous section, with
some minor modifications. By virtue of Sibuya's theorem,
equation (\ref{donthuc}) possesses a unique solution
$\Phi_0(X,\lambda)$, which is an entire function in both $X$ and
$\lambda$. For each fixed $\lambda \in \C$, this solution and its
derivative with respect to $X$ satisfy the following asymptotic
estimates
\begin{equation}\label{caseI}
    \begin{array}{ll}
       &  \Phi_0(X,\lambda) \simeq X^{-m/4} e^{-\frac{2}{m+2}X^{\frac{m+2}{2}}} \left[
       1+O(X^{-1/2})\right]\\\vspace{2mm}
      {\rm and} & \Phi'_0(X,\lambda) \simeq X^{m/4} e^{-\frac{2}{m+2}X^{\frac{m+2}{2}}} \left[
  -1+O(X^{-1/2})\right] \\
    \end{array}
    \qquad {\rm as}\quad
      S_0 \ni X \rightarrow \infty
\end{equation}
We should notice that, both $\Phi_0$ and $\Phi'_0$ vanish
exponentially at infinity in the sector $S_0$. Besides, as an
entire function of $X$, $\Phi_0(X,\lambda)$ can be expanded in
powers of $X$,
\begin{equation}\label{sohangtudo}
\Phi_0(X,\lambda) =
f_0(\lambda)+\sum_{j=1}^{\infty}f_j(\lambda)X^j,
\end{equation}
where $f_j(\lambda)$ are also entire functions of $\lambda$ and
satisfy the symmetry property
$$\overline{f_j(\lambda)}= f_j(\overline{\lambda})\, ,\quad j = 0,1,2,\ldots$$
The companion solutions $\Phi_{\pm 1}$ of $\Phi_0$ now read as
follows:
\begin{equation}\label{}
    \begin{array}{l}
        \Phi_1(X,\lambda) = \omega^{-m/4}\Phi_0(\omega^{-1}X,\omega^{2}\lambda)=\omega^{-m/4}
       \big\{f_0(\omega^2\lambda)+\sum_{j=1}^{\infty}f_j(\omega^2\lambda)(\omega^{-1}X)^j\big\}
       \\\vspace{3mm}
      \Phi_{-1}(X,\lambda) = \omega^{m/4}\Phi_0(\omega X,\omega^{-2}\lambda)=\omega^{m/4}
       \big\{f_0(\omega^{-2}\lambda)+\sum_{j=1}^{\infty}f_j(\omega^{-2}\lambda)(\omega X)^j\big\} \\
    \end{array}
\end{equation}
Substituting these expressions into (\ref{ourcasedx}) and
specializing that equality with $X= 0$ lead the following relation
between the Stokes multiplier $C(\lambda)$ and the entire function
$f_0(\lambda)$:
\begin{equation}\label{hethuccoban}
    C(\lambda)f_0(\lambda) = \omega^{m/4}f_0(\omega^{-2}\lambda)-\omega^{-m/4}f_0(\omega^{2}\lambda).
\end{equation}
This intimate relation can serve to study the zeros of
$C(\lambda)$ through those of $f_0(\lambda)$. Before going
further, we recall a consequence from Sibuya's asymptotic studies
 on these functions (\cite{Sibuya}, ch. $4,5$).
\begin{Pro}\label{capsohangtudo}
The orders of entire functions $f_0(\lambda)$ and  $C(\lambda)$
are both equal to  $\frac{1}{2}+\frac{1}{m}\,$.
\end{Pro}
\begin{Rem}\label{}
    For $m \geq 3$, the order is not an integer.
    Therefore, $f_0(\lambda)$ and $C(\lambda)$ must have infinitely many
    zeros, whose accumulation point can only be infinity.
\end{Rem}

The following assertion locates the zeros of $f_0(\lambda)$ more
concretely.
\begin{Pro}\label{zerosf0}
All the zeros of $f_0(\lambda)$ are negative real numbers.
Moreover,
$$f_0(\lambda) > 0\, , \qquad \forall \lambda \geq 0$$
\end{Pro}
Before proving the proposition, we need to remind a useful transform for a given
second-order linear differential equation in the complex plane.
Let $w(z)$ be a solution of the following equation,
\begin{equation}
  \label{ptvpcaphai}
  w''(z) - f(z)w(z) = 0.
\end{equation}
Then for any $z_1, z_2 \in \C$, we have the following identity
by multiplying (\ref{ptvpcaphai}) by $\overline{w(z)}$ and
integrating it from $z_1$ to $z_2$:
\begin{equation}\label{Green}
\overline{w(z)}w'(z)\Big|_{z_1}^{z_2} := \overline{w(z_2)}w'(z_2)
- \overline{w(z_1)}w'(z_1) =
  \int_{z_1}^{z_2}|w'(z)|^2d\overline{z} +
\int_{z_1}^{z_2}f(z)|w(z)|^2 dz
\end{equation}
This identity, which is known as the \emph{Green's transform} of
(\ref{ptvpcaphai}), is true provided that the integrals in the
right-hand side of (\ref{Green}) make sense.

In case where the integral path is a segment $[z_1, z_2] = \{z(t)=z_1 +
te^{i\theta}/\, t \in [0,r]\}$, with $\theta = \arg(z_2-z_1)$ and
$r = |z_2-z_1|$, then (\ref{Green}) turns into
\begin{equation}
  \label{Greentrendoan}
  \overline{w(z)}w'(z)\Big|_{z_1}^{z_2}
= e^{-i\theta} \int_0^r|w'(z(t))|^2dt +
e^{i\theta}\int_0^r|w(z(t))|^2f(z(t))dt.
\end{equation}
Besides, this equality holds true as  $r\rightarrow +\infty$,
provided all limits exist.

\vspace{4mm}\noindent\emph{Proof of Proposition \ref{zerosf0}.}
    Let $\lambda^* = \alpha+i\beta$ be a zero of $f_0(\lambda) = \Phi_0(0,\lambda)$.
    By applying the Green's transform (\ref{Greentrendoan}) on the interval $[0,X] \subset \R$
    to the solution $\Phi_0$ of
    (\ref{donthuc}) , we obtain
    \begin{equation}\label{apdungGreen}
    \overline{\Phi_0(t,\lambda)}\Phi'_0(t,\lambda)\Big|_{0}^{X}
= \int_0^{X}|\Phi'_0(t,\lambda)|^2dt +
\int_0^{X}(t^m+\lambda)|\Phi_0(t,\lambda)|^2dt.
\end{equation}
Substituting $\lambda = \lambda^*$ into (\ref{apdungGreen}) and
letting $X\rightarrow +\infty$ yield
$$0 = \int_0^{+\infty}|\Phi'_0(t,\lambda^*)|^2dt +
\int_0^{+\infty}(t^m+ \alpha+i\beta)|\Phi_0(t,\lambda^*)|^2dt
$$
By separating the real and imaginary parts, we obtain
$$\beta = 0 \quad {\rm and }\quad \alpha =
 -\frac{\int_0^{+\infty}|\Phi'_0(t,\lambda^*)|^2dt+\int_0^{+\infty}t^m|\Phi_0(t,\lambda^*)|^2dt}
 {\int_0^{+\infty}|\Phi_0(t,\lambda^*)|^2dt} < 0$$
For the rest of the proposition, we consider a fixed $\lambda \geq
0$ and treat $X$ as a real variable. Hence,  equality
(\ref{doixung}) implies that $\Phi_0(X,\lambda)$ is a real-valued function of
$X\in \R$, and so is also its derivative $\Phi'_0(X,\lambda)$.
By taking the derivative of (\ref{apdungGreen}) with
respect to $X$, we get
$$\big(\overline{\Phi_0(X,\lambda)}\Phi'_0(X,\lambda)\big)'_X = |\Phi'_0(X,\lambda)|^2 +
(X^m+\lambda)|\Phi_0(X,\lambda)|^2 .
$$
Since $\lambda \geq 0$, the right-hand side of this equality is
strictly positive on $[0,+\infty)$. It follows that the real
function $\Phi_0(X,\lambda)\Phi'_0(X,\lambda)$ is strictly
increasing on $[0,+\infty)$. Moreover, we deduce from
(\ref{caseI}) that
$$\lim_{X\rightarrow+\infty}\Phi_0(X,\lambda)\Phi'_0(X,\lambda)=0$$
Combining these facts implies that both $\Phi_0(X,\lambda)$ and
$\Phi'_0(X,\lambda)$ never vanish on $[0,+\infty)$.

In particular, regarding its asymptotic asymptotic behavior in
(\ref{caseI}), we can conclude that $\Phi_0(X,\lambda) > 0$ for
all $X\geq 0$. Putting $X=0$ completes the proof.\dpcm

Concerning the zeros of $C(\lambda)$, we have the following.
\begin{Thm}\label{tinhdon1}
All  zeros of $C(\lambda)$ are real, positive and simple.
\end{Thm}
We should remind that the reality of all these zeros, in
connection with $\mathcal{PT}$-symmetric quantum mechanics, has
been studied in various ways by many authors
\cite{BenBoe98,DP98-1,DT00,Dorey,Shin002,Tai2}. The simpleness of all
but a finite number of zeros has been indicated in \cite{Sibuya} (ch.
$6$) by using some asymptotic estimates for the large zeros.

Next, we shall provide a rigorous proof for the simpleness of all
zeros using ideas from the proof of Laguerre's theorem
\cite{Boas}, after simply justifying the reality and positivity
 in our special case by the same way as in \cite{Shin002}.

\vspace{4mm}\noindent\emph{Proof of Theorem \ref{tinhdon1}.}  Let
$\{\lambda_n\}_{n=\overline{0,\,\infty}}$ be zeros of the entire
function $f_0(\lambda)$. Applying  Hadamard's factorization
theorem (\cite{Boas}) to $f_0(\lambda)$, whose order is smaller than $1$ for
$m>2$, we have, for all $\lambda \in \mathbb{C}$,
\begin{equation}\label{Hadamard}
    f_0(\lambda) = A\prod_{n=0}^{\infty}\Big(1-\frac{\lambda}{\lambda_n}\Big)
\end{equation}
where $A= f_0(0) > 0$.

Suppose that $\lambda^*$ is a zero of $C(\lambda)$. We can deduce
from (\ref{hethuccoban}) that
\begin{equation}\label{lalala}
A\omega^{m/4}\prod_{n=0}^{\infty}\Big(1-\frac{\omega^{-2}\lambda^*}{\lambda_n}\Big)
=A\omega^{-m/4}\prod_{n=0}^{\infty}\Big(1-\frac{\omega^{2}\lambda^*}{\lambda_n}\Big)
\end{equation}
By virtue of Proposition \ref{zerosf0},  all the zeros $\lambda_n $
of $f_0(\lambda)$ are negative. Combining this with the fact that
$0<|\arg(\omega^{\pm 2})| <\pi$, we can conclude that
$f_0(\omega^{\pm 2}\lambda^*)$ cannot be simultaneously vanishing.
Therefore, both of sides of (\ref{lalala}) are never vanishing.

By taking the absolute value of both sides and remarking $|z| =
|\overline{z}|$, we obtain
\begin{equation}\label{tamthoi}
\prod_{n=0}^{\infty}\Big|\frac{\omega^{2}\lambda_n-\lambda^*}{\omega^{2}\lambda_n-\overline{\lambda^*}}\Big|=1
\end{equation}
Since $\lambda_n <0$, we have ${\rm Im}(\lambda_n\omega^{2}) < 0,
\, \forall n$. These inequalities  imply that, unless $\lambda^* =
\overline{\lambda^*}$, the factors in the left hand side of
(\ref{tamthoi}) are all either strictly greater or less than $1$.
Consequently, the truth of the equality (\ref{tamthoi}) requires
that $\lambda^* \in \R$.

To verify that $\lambda^* > 0$, we use the Green's transform again.
Applying (\ref{Greentrendoan}) to the solution
$\Phi_1(X,\lambda^*)$ of the equation (\ref{donthuc}) on the ray
$[0,\omega\infty)$,  we obtain
\begin{equation}
  \label{tam2}
  \overline{\Phi_1(X,\lambda^*)}\Phi'_1(X,\lambda^*)\Big|_{0}^{\omega\infty}
= \omega^{-1}\int_0^{\infty}|\Phi'_1(\omega t,\lambda^*)|^2dt +
\omega\int_0^{\infty}|\Phi_1(\omega t,\lambda^*)|^2((\omega t)^m
+\lambda^*)dt
\end{equation}
Note that, by definition, $\Phi_1(X,\lambda^*)$ and its derivative
$\Phi'_1(X,\lambda^*)$ are exponentially vanishing at
$\omega\infty $. On the other hand, since $\lambda^*$ is a (real)
zero of $C(\lambda)$, we can deduce from (\ref{ourcasedx}) and
(\ref{dxcuanghiem}) that
$$\Phi_1(X,\lambda^*) = \Phi_{-1}(X,\lambda^*) =
\overline{\Phi_1(\overline{X},\overline{\lambda^*})}
$$ This implies that the left-hand side of (\ref{tam2}),
which now becomes
$-\overline{\Phi_1(0,\lambda^*)}\Phi'_1(0,\lambda^*)$, is purely
real.

Separating the imaginary part in (\ref{tam2}), where
$\omega=e^{i\theta}$ and $\theta=\frac{2\pi}{m+2}$, we obtain
$$0= -\sin\theta\int_0^{\infty}|\Phi'_1(\omega t,\lambda^*)|^2dt
-\sin\theta\int_0^{\infty}|\Phi_1(\omega t,\lambda^*)|^2t^m dt +
\lambda^*\sin\theta\int_0^{\infty}|\Phi_1(\omega
t,\lambda^*)|^2dt$$ This equality indicates the positivity of
$\lambda^*$.

To finish the demonstration, we have to show that $C'(\lambda)
\neq 0$ if $C(\lambda) =0$. Since all zeros of $C(\lambda)$ are
real according to the above proofs, it is sufficient to consider
$\lambda$ as a real variable.

We now turn to (\ref{hethuccoban}). By setting
$$g(\lambda) = \omega^{m/4}f_0(\omega^{-2}\lambda)$$
and regarding it as a complex-valued function of the real variable
$\lambda$, we can write
$$g(\lambda) = g_1(\lambda) + ig_2(\lambda)$$
where $g_{1,2}(\lambda)$ are real  and differentiable functions on
$\R$.

On account of the relation (\ref{hethuccoban}), we see that a
double zero $\lambda^*\in\R$ of $C(\lambda)$ must satisfy
$$g(\lambda^*) = \overline{g(\lambda^*)}\quad {\rm and}\quad
g'(\lambda^*) = \overline{g'(\lambda^*)}. $$
This is equivalent to requiring that
$$
g_2(\lambda^*) =0\quad {\rm and}\quad g'_2(\lambda^*) =0,
$$
that is (since we know that $g(\lambda^*) \neq 0$),
\begin{equation}\label{tam4}
\frac{g'(\lambda^*)}{g(\lambda^*)}\in \R.
\end{equation}
But, as indicated below, this is impossible.

Indeed, using (\ref{Hadamard}), we can factorize $g(\lambda)$ as
follows, for $\lambda \in \mathbb{R}$,
$$g(\lambda) = A\omega^{m/4}\prod_{n=0}^{\infty}\Big(1-\frac{\omega^{-2}\lambda}{\lambda_n}\Big)$$
This factorization also indicates that $g(\lambda) \neq 0$ on
$\R$. Taking the logarithmic derivative on both sides of the equality
yields
\begin{equation}\label{tam3}
\frac{g'(\lambda)}{g(\lambda)} =
\sum_{n=0}^{\infty}\frac{1}{\lambda-\lambda_n\omega^2}.
\end{equation}
Moreover, we have  $\lambda_n <0$. Hence, for all $\lambda\in\R$,
the imaginary part of the right-hand side of (\ref{tam3}) is
$$\sin2\theta\sum_{n=0}^{\infty}\frac{\lambda_n}{(\lambda-\lambda_n\cos2\theta)^2 +(\lambda_n\sin2\theta)^2}$$
and it is never vanishing; against (\ref{tam4}). This
completes the proof.\dpcm

\begin{Rem}\label{rem1}
By following the above proof, we can conclude that all the zeros of
  $f_0(\lambda)$ are also simple, in connection with
  Proposition~\ref{zerosf0} affirming the reality of these zeros.
\end{Rem}
    As a consequence of  this theorem, the real-valued function $-iC'(\lambda)$ changes its sign alternately
    at zeros of $C(\lambda)$. This matter has been mentioned in \cite{Tai3}
    as an attempt to justify the indefiniteness of
    $\mathcal{PT}$-pseudo-norm in $\mathcal{PT}$-symmetric quantum
    mechanics.

\subsection{General case}
In what follows, we proceed with the study of the zeros of the Stokes
multiplier in the case where $a_k$ are not simultaneously equal to
zero. Since we shall use again the previous arguments, a brief recall
of the notions should be done.

For our goal, we consider the following equation:
\begin{equation}\label{tongquat}
-\Phi''(X) + (X^{m} +a_{1}X^{m-1}+\cdots +a_{m-1}X
+\lambda)\Phi(X)=0
\end{equation}
We start with one of its solutions $\Phi_0(X,a,\lambda)$, whose
existence and asymptotic behavior have been settled in the
Sibuya's theorem \ref{ketquaSibuya}. As an entire function of $X$,
$\Phi_0$ can be written in the form
\begin{equation}\label{tam5}
\Phi_0(X,a,\lambda) = f_0(a,\lambda) +
\sum_{j=1}^{\infty}f_j(a,\lambda)X^j
\end{equation}
where $f_0(a,\lambda)$ is an entire function in both $a =
(a_1,\ldots,a_{m-1})$ and $\lambda$. In particular, considered as
an entire function in $\lambda$, $f_0(a,\lambda)$ is of order
 $\frac{1}{2}+\frac{1}{m}$ (uniformly in $a$ for $a$ in a compact
 set).
The companions solutions $\Phi_{\pm 1}$
are defined by
\begin{equation}\label{}
    \begin{array}{l}
        \Phi_1(X,a,\lambda) = \omega^{-m/2-r_m(a)}
        \Phi_0(\omega^{-1}X,\omega_{-1}(a),\omega^{-m}\lambda)\\\vspace{3mm}
      \Phi_{-1}(X,a,\lambda) = \omega^{m/2+r_m(a)}
        \Phi_0(\omega X,\omega_{1}(a),\omega^{m}\lambda) \\
    \end{array}
\end{equation}
The relation among these three solutions is realized  by the Stokes multiplier $C(a,\lambda)$:
$$\Phi_{-1}(X,a,\lambda) = C(a,\lambda)\Phi_{0}(X,a,\lambda) +\Phi_{1}(X,a,\lambda).$$
Putting $X=0$ in this equality, together with (\ref{tam5}), we
obtain
\begin{equation}\label{hethuc2}
C(a,\lambda)f_{0}(a,\lambda) =
\omega^{m/2+r_m(a)}f_{0}(\omega_{1}(a),\omega^m\lambda) -
\omega^{-m/2-r_m(a)}f_{0}(\omega_{-1}(a),\omega^{-m}\lambda).
\end{equation}
For a fixed $a = (a_1,\dots,a_{m-1})\in\R^{m-1}$, we are examining
the zeros of the first term\footnote{The other one is nothing but
$\overline{g_a(\overline{\lambda})}$\,.}  in the right-hand side
of (\ref{hethuc2})
$$ g_a(\lambda) := \omega^{m/2+r_m(a)}f_{0}(\omega_{1}(a),\omega^m\lambda).$$
Note that $g_a(\lambda)$ is also an entire function of order
$\frac{1}{2}+\frac{1}{m}$, so for $m>2$, the order is
non-integral. It follows that $g_a(\lambda)$ must have an infinite
number of zeros.

Let $\lambda_n = \lambda_n(a)$, $n \in \N$, be the zeros of
$g_a(\lambda)$. Then $\Phi_0(0,\omega_{1}(a),\omega^m\lambda_n) =
0$, where $ Y(X) := \Phi_0(X,\omega_{1}(a),\omega^m\lambda_n)$
verifies the following equation
\begin{equation}\label{tam6}
    -Y''(X) + (X^{m} +a_{1}\omega^1X^{m-1}+\cdots +a_{m-1}\omega^{m-1}X
+\lambda_n\omega^m)Y(X)=0
\end{equation}
We now briefly discuss a significant result of Shin on the reality
of the zeros of $C(a,\lambda)$ in \cite{Shin002}, where the author
give a suitable hypothesis on the real coefficients $a_k$ under which
all zeros of $C(a,\lambda)$ are showed to be real and positive.
That hypothesis reads
$$\mbox{\sf There exists an index}\ j,\, 1\leq j \leq m/2 \ \mbox{\sf such that }\,(j-k)a_k\geq 0
\ \mbox{\sf for all } k\,.  \eqno{(H)}
$$
One of the crucial steps in his proof is to effectively apply the
Green's transform  to (\ref{tam6}) on a suitably chosen ray,  so
that all the imaginary parts of $\lambda_n$ are non-positive.
 For our present purpose, we require that all ${\rm Im}\lambda_n$ be
 strictly negative. This fact  can actually be derived from the
 hypothesis  $(H)$ thanks to Shin's proof itself
 in the cited article, the special case when $m=4$ and $j=2$ being apart.
To overcome it, we only need to add  a negligible supplement on
$(H)$ as follows
$$\mbox{\sf If}\ (H)\ \mbox{\sf occurs for}\ m=4\ {\sf and}\ j=2,\ {\sf then}\ a_2 \leq 0\, .\eqno{(s)}
$$
Indeed, for $m=4$ the hypothesis $(H)$ and its supplement $(s)$
imply that $a_1 \geq 0$ and $a_{2,3} \leq 0$.  Applying
(\ref{Greentrendoan}) to (\ref{tam6}) on the ray
$[0,e^{i\theta}\infty)$, where $\omega = e^{\frac{i\pi}{3}}$ and
$|\theta| < \frac{\pi}{6}$, we obtain
$$\begin{array}{lll}
    0 & = & \displaystyle  e^{-i\theta}\int_0^{+\infty}|Y'|^2dt +
    \\ \vspace{2mm}
      &  & \displaystyle + e^{i\theta}\int_0^{+\infty}|Y|^2(e^{i4\theta}t^4 + a_1\omega
e^{i3\theta}t^3
 + a_2\omega^2
e^{i2\theta}t^2+a_3\omega^3e^{i\theta}t + \lambda_n\omega^4)dt
\end{array}
$$
Taking the imaginary part in this equality, after multiplying it
by $e^{i(-\theta+2\pi/3)}$, we get
\begin{equation}\label{bs}
\begin{array}{lll}
  0 &= & b_0\sin(2\pi/3-2\theta)+b_1\sin(2\pi/3+4\theta) + b_2a_1\sin(\pi +3\theta) + \\
  \vspace{3mm}
   & &  + b_3a_2\sin(4\pi/3+2\theta) + b_4a_3\sin(5\pi/3+\theta) +b_5{\rm Im}\lambda_n \\
\end{array}
\end{equation}
where the constants  $b_j >0$ stand for the values of the
integrals.

Following the above conditions, we can deduce that\footnote{In
fact, it holds  for any real $a_1$.}  ${\rm Im}\lambda_n < 0$ by
letting $\theta =0$.

 We now come to the following theorem, which is
more general than theorem \ref{tinhdon1}.
\begin{Thm}\label{dinhly2}
    For a fixed $a = (a_1,\ldots,a_{m-1}) \in \R^{m-1}$ which satisfies the
    hypothesis $(H)$ and its supplement $(s)$, all the zeros of the  Stokes
    multiplier $C(a,\lambda)$ are real, positive and simple.
\end{Thm}
\begin{proof}\label{}
    The reality and positivity have been already proved by Shin \cite{Shin002}.
The simpleness can be handled in the same manner as for the proof of
theorem \ref{tinhdon1}, thus allowing us to be sketchy.

Let $\lambda^* = \lambda^*(a)$ be a  double zero of
$C(a,\lambda)$. Because of the reality of the zeros, it is sufficient
to consider $C(a,\lambda)$ as a function of $\lambda\in \R$. We
then deduce from (\ref{hethuc2}) that
\begin{equation}\label{tam8}
 g_a(\lambda^*) =
\overline{g_a(\lambda^*)}\quad \mbox{and}\quad
(g_a)'_{\lambda}(\lambda^*) =
\overline{(g_a)'_{\lambda}(\lambda^*)}
\end{equation}
Note that $g_a(\lambda)$ is an entire function in $\lambda$ of
order $\frac{1}{2}+\frac{1}{m} < 1$, whose zeros $\lambda_n$ now
satisfy ${\rm Im}\lambda_n <0$. By the Hadamard's factorization
theorem, we have
$$g_a(\lambda) = g_a(0)\prod_{n=0}^{\infty}\Big(1-\frac{\lambda}{\lambda_n}\Big).$$
Obviously, $g_a(\lambda) \neq 0$ on $\R$, so taking the
logarithmic derivative this identity yields
$$\frac{(g_a)'_{\lambda}(\lambda)}{g_a(\lambda)} = \sum_{n=0}^{\infty}\frac{1}{\lambda-\lambda_n}.$$
The imaginary part of the right hand member is
$$\sum_{n=0}^{\infty}\frac{{\rm Im}\lambda_n}{(\lambda-{\rm Re}\lambda_n)^2+({\rm Im}\lambda_n)^2} <0.$$
Thus, (\ref{tam8}) cannot be realized at any $\lambda^*\in \R$.
\end{proof}

\section{Conclusion}

In this paper, we have shown that, under an appropriate hypothesis of
signs on the real coefficients $a_1,\ldots,a_{m-1}$ of equation
(\ref{EquatX}), all the zeros of the Stokes multiplier $C(a)$ are
simple. On the one hand, our proof is partly based on the results of
\cite{Shin002}, in particular the hypothesis that we use to ensure the
reality of these zeros is a sufficient condition  first established by
Shin.
On the other hand, our main  remaining arguments, essentially the
 Green's transform and the Hadamard's factorization theorem, are
 certainly natural in this context, being already used in various
 papers. In particular, whenever no
parameter is concerned, our reasonings are quite simple.

Since our main theorem \ref{dinhly2} makes use of a sufficient for the
zeros to be real, we do not give any information about neither
the simpleness of the zeros for the cases where some of them are
complex, nor conditions on the parameter $a$ for the existence of
multiple (real or complex) zeros. As an illustration of the case exhibiting double real
zeros, we suggest a common paper with Delabaere \cite{DT00}, where
the energy spectrum of the Hamiltonian $H = p^2+i(q^3 + \alpha q)$
acting on $L^2(\R)$ was studied by semiclassical analysis. The
differential equation associated with this Hamiltonian reads
$-\Phi''(X) + (X^3 -\alpha X+E)\Phi(X) =0$, up to a rotation (see
also \cite{Tai2}). It has been shown that in this case some pairs
of real zeros $E_n(\alpha)$ of the Stokes multiplier may coalesce
before turning into complex conjugate at certain critical values
of $\alpha_{crit} < 0$. These degenerate values
$E_n(\alpha_{crit})$ are nothing but common zeros of the Stokes
multiplier and its derivative with respect to $E$ (see Fig.1 in
\cite{DT00}). It seems to us that the question of characterizing these critical values
of $\alpha$ is certainly an interesting but quite challenging problem.

From the mathematical viewpoint, this question is of course
related to the fact that in general for $m \geq 3$, neither
special functions solution of (\ref{EquatX}) are known, nor their
related Stokes multipliers (when $m=2$  these multipliers can be
explicitly expressed in term of the Gamma function, while for
$m=1$ they are constants). In this way, the present paper can be
thought of as an attempt for exploring some hidden special
functions in relation  to their Stokes multipliers. But, as
already said in the introduction, the main motivation for this
paper was to add a new (small) stone for the mathematical
foundation of   $\mathcal{PT}$-symmetric models, in addition to
our results in \cite{Tai3}.

\vspace*{5mm}

\subsection*{Acknowledgments.} This work was
supported by the International Centre for Theoretical Physics in
the framework of Post-doctoral Fellowship.


\end{document}